\documentclass[11pt, reqno]{amsart}%
\usepackage{amsmath, amstext, amsbsy, amssymb, amscd}
\usepackage{amsmath}
\usepackage{amsxtra}
\usepackage{amscd}
\usepackage{amsthm}
\usepackage{amsfonts}
\usepackage{amssymb}
\usepackage{eucal}
\usepackage{color}
\usepackage{graphicx}%
\newtheorem{theorem}{Theorem}[section]
\theoremstyle{plain}

\newtheorem{lemma}[theorem]{Lemma}
\newtheorem{proposition}[theorem]{Proposition}
\newtheorem{corollary}[theorem]{Corollary}

\theoremstyle{definition}
\newtheorem{definition}[theorem]{Definition}

\theoremstyle{remark}
\newtheorem{remark}[theorem]{Remark}

\definecolor{A}{rgb}{.75,1,.75}

\numberwithin{equation}{section}

\newcommand{\C}{\mathbb C}
\newcommand{\Z}{\mathbb Z}
\newcommand{\N}{\mathbb N}
\newcommand{\Cl}{\mathcal C}
\newcommand{\mh}{\mathfrak{h}}          

\newcommand{\be}{\beta}

\newcommand{\al}{\alpha}

\newcommand{\td}{\tilde}
\newcommand{\disp}{\displaystyle}

\newcommand{\wtd}{\widetilde}

\newcommand{\aHC}{\mathfrak{H}^{\mathfrak c}} 
\newcommand{\daHC}{\ddot{\mathfrak{H}}^{\mathfrak c}}
\newcommand{\sdaH}{\ddot{\mathfrak{H}}^-}
\newcommand{\saH}{\mathfrak{H}^-} 
\newcommand{\sDaH}{\ddot{H}^{\mathfrak{-}}}
\newcommand{\dDaH}{\ddot{H}} 
\newcommand{\dDaHC}{\ddot{H}^{\mathfrak c}}

{\vskip-\lastskip\medskip
  \noindent
  {\em #1.}\enspace
  }%
{\qed\par\medskip
  }

\begin{document}

\title[The classical trigonometric spin DAHA]
{Hecke-Clifford algebras and spin Hecke algebras~III: the trigonometric type}

\author[Ta Khongsap]{Ta Khongsap}
\address{Department of Math., University of Virginia,
Charlottesville, VA 22904} \email{tk7p@virginia.edu}

\begin{abstract}
The notion of trigonometric spin double affine Hecke algebras (tsDaHa)
and trigonometric double affine Hecke-Clifford algebras (tDaHCa)
associated to classical Weyl groups are introduced. The PBW basis property is established. An algebra isomorphism relating tDaHCa to tsDaHa is obtained.
\end{abstract}

\maketitle

\section{Introduction}\label{intro}

\subsection{}
For an irreducible finite Weyl group $W$ associated to a root system $R$, there are corresponding affine Weyl group $W^a$ and extended affine Weyl group $W^e$ attached with an affine root system $\wtd{R}$. There is an interesting family of algebras attached to $W^a$ or $W^e$, namely the trigonometric double affine Hecke algebra (tDAHA), $\dDaH_{t,c}$, where $t,c$ are certain parameters. The tDAHA is not only a degeneration of the double affine Hecke algebra (DAHA), but also an extension of the degenerate affine Hecke algebra (AHA) as defined by Lusztig in \cite{Lu}. The tDAHA was introduced and studied by Cherednik with application to harmonic analysis and Macdonald polynomials, see \cite{Ch1, Ch2}. Furthermore, when one specializes $t=0$, the algebra $\dDaH_{t=0,c}$ has a large center, in particular, $\dDaH_{t=0,c}$ is a finite module over its center, and it has interesting connections with algebraic geometry (see \cite{Ob}).

In 1911, I. Schur developed a theory of spin (projective) representations of the symmetric group $S_n$. We note that studying the spin representations of $S_n$ is equivalent to study the representations of the \textit{spin} symmetric group algebra $\C S_n^-$. For the symmetric group $S_n$, there is a standard procedure to construct the associated Hecke algebras. Wang \cite{W} has raised the question of whether or not a notion of Hecke algebras associated to the \textit{spin} symmetric group algebra $\C S_n^-$ exists, and has provided a natural construction of the trigonometric DAHA associated to the algebra $\C S_n^-$, denoted by $\sdaH_{tr}$. Moreover, in  \cite{W,KW1,KW2}, we have developed a theory of spin Hecke algebras of (degenerate) affine and rational double affine type (also cf. \cite{Naz}).

\subsection{}
This paper is a sequel to \cite{W,KW1,KW2}. The main goal here is to construct two classes of algebras which are closely related to the tDAHA associated to each classical finite Weyl group $W$ called the \textit{trigonometric double affine Hecke-Clifford algebra} (tDaHCa), $\dDaHC_W$, and the \textit{trigonometric spin double affine Hecke algebra} (tsDaHa), $\sDaH_W$. We establish a few basic properties of the algebras $\dDaHC_W$ and $\sDaH_W$, including their PBW basis properties. In addition, we prove that for type $A$ case, the algebra $\dDaHC_W$ (respectively, $\sDaH_W$) is isomorphic to $\daHC_{tr}$ (respectively, to $\sdaH_{tr}$) introduced in \cite{W} in a very different presentation.

\subsection{}
In Section \ref{spin:extend Weyl}, we recall the definition of an extended affine Weyl group $W^e$ for each classical type, and then formulate the corresponding \textit{spin} extended affine Weyl group algebra $\C W^{e-}$. We first start with a finite Weyl group $W$, and then consider a double cover $\wtd{W}$ of $W$ associated to a distinguished $2$-cocycle as in \cite{KW1,Mo}:
\begin{eqnarray*} \label{ext}
1 \longrightarrow \Z_2 \longrightarrow \wtd{W} \longrightarrow W
\longrightarrow 1.
\end{eqnarray*}
We then define a covering of the extended affine Weyl group $\wtd{W}^e$ of $W^e$. As a result, we obtain the \textit{spin} extended affine Weyl group algebra $\C W^{e-}$ which is defined to be a certain quotient algebra of $\C \wtd{W}^e$.

In Section \ref{Clifford:Extended Weyl}, recalling that the Weyl group $W$ acts as automorphisms on the Clifford algebra $\Cl_n$, we then extend  to an action of $W^e$ on $\Cl_n$. We establish a (super)algebra isomorphism
\[
    \Phi: \Cl_n \rtimes \C W^e \stackrel{\simeq}{\longrightarrow} \Cl_n \otimes \C W^{e-}
\] which is an extended affine analogue to the isomorphism $~\Phi: \Cl_n \rtimes \C W \stackrel{\simeq}{\longrightarrow} \Cl_n \otimes \C W^{-}$ in \cite{KW1} which goes back to Sergeev \cite{Ser} and Yamaguchi \cite{Yam} for type A.

In Section \ref{d DaHC} we introduce the tDaHCa $\dDaHC_W$. We establish the PBW basis properties for the algebras $\dDaHC_W$:
$$
\dDaHC_W \cong \C [\mh^*] \otimes \Cl_n \otimes \C W^e
$$
where $\C [\mh^*]$ denotes the polynomial algebra.

In Section \ref{sec:spin}, we introduce the tsDaHa $\sDaH_W$. The tsDaHa are not only the counterpart of the algebras $\dDaHC_W$, but also the generalizations of the (degenerate) spin affine Hecke algebras associated to the \textit{spin} group algebras $\C W^-$ in \cite{W,KW1}. Denote by $\Cl [\mh^*]$ a noncommutative skew-polynomial algebra, we establish the PBW basis properties for $\sDaH_W$:
$$
\sDaH_W \cong \Cl [\mh^*] \otimes \C W^{e-}.
$$
In addition, we establish a (super)algebra isomorphism:
$$
\Phi: \dDaHC_W \stackrel{\simeq}{\longrightarrow} \Cl_n \otimes
\sDaH_W
$$ which extends the isomorphism $\Phi: \Cl_n \rtimes \C W^e \stackrel{\simeq}{\longrightarrow} \Cl_n \otimes \C W^{e-}$.

\subsection{}
The constructions of our algebras are canonical, but we have to do it case-by-case since the algebras in a way rely on a choice of orthonormal basis of $\mh$. However, we hope that the detailed presentations on extended affine Weyl groups could be helpful to the reader.

{\bf Acknowledgements.} I am thankful for the Semester Dissertation Fellowship which allowed me to concentrate on my research during the semester of Spring 2008. I am deeply grateful to my advisor W. Wang for many useful suggestions during this project.
%
%
%
%
%
%
%
\section{(Spin) extended affine Weyl Groups} \label{spin:extend Weyl}

\subsection{The Weyl group $W$}
Let $W$ be an (irreducible) finite Weyl group of classical type with the following presentation:
\begin{eqnarray} \label{eq:weyl}
\langle s_1,\ldots,s_n | (s_is_j)^{m_{ij}} = 1,\ m_{i i} = 1,
 \ m_{i j} = m_{j i} \in \Z_{\geq 2}, \text{for } i
 \neq j \rangle.
\end{eqnarray}
The integers $m_{i j}$ are specified by the Coxeter-Dynkin diagrams
whose vertices correspond to the generators of $W$ below. By
convention, we only mark the edge connecting $i,j$ with $m_{ij} \ge
4$. We have $m_{ij}=3$ for $i \neq j$ connected by an unmarked edge,
and $m_{ij}=2$ if $i,j$ are not connected by an edge.

 \begin{equation*}
 \begin{picture}(150,45) 
 \put(-99,18){$A_{n}$}
 \put(-30,20){$\circ$}
 \put(-23,23){\line(1,0){32}}
 \put(10,20){$\circ$}
 \put(17,23){\line(1,0){23}}
 \put(41,22){ \dots }
 \put(64,23){\line(1,0){18}}
 \put(82,20){$\circ$}
 \put(89,23){\line(1,0){32}}
 \put(122,20){$\circ$}
 \put(-30,9){$1$}
 \put(10,9){$2$}
 \put(74,9){${n-1}$}
 \put(122,9){${n}$}
 \end{picture}
 \end{equation*}
 %
 \begin{equation*}
 \begin{picture}(150,55) 
 \put(-99,18){$B_{n}(n\ge 2)$}
 \put(-30,20){$\circ$}
 \put(-23,23){\line(1,0){32}}
 \put(10,20){$\circ$}
 \put(17,23){\line(1,0){23}}
 \put(41,22){ \dots }
 \put(64,23){\line(1,0){18}}
 \put(82,20){$\circ$}
 \put(89,23){\line(1,0){32}}
 \put(122,20){$\circ$}
 \put(-30,10){$1$}
 \put(10,10){$2$}
 \put(74,10){${n-1}$}
 \put(122,10){${n}$}
 %
 \put(102,24){$4$}
 \end{picture}
 \end{equation*}
%
%
 \begin{equation*}
 \begin{picture}(150,75) 
 \put(-99,28){$D_{n} (n \ge 4)$}
 \put(-30,30){$\circ$}
 \put(-23,33){\line(1,0){32}}
 \put(10,30){$\circ$}
 \put(17,33){\line(1,0){15}}
 \put(35,30){$\cdots$}
 \put(52,33){\line(1,0){15}}
 \put(68,30){$\circ$ }
 \put(75,33){\line(1,0){32}}
 \put(108,30){$\circ$}
 \put(113,36){\line(1,1){25}}
 \put(138,61){$\circ$}
 \put(113,29){\line(1,-1){25}}
 \put(138,-1){$\circ$}
 \put(-29,20){$1$}
 \put(10,20){$2$}
 \put(60,20){$n-3$}
 \put(117,30){$n-2$}
 \put(145,0){$n-1$}
 \put(145,60){$n$}
 %
 \end{picture}
 \end{equation*}
\vspace{.5cm}
Denote by $W_{A_{n-1}}$ (respectively, $W_{B_n}$ and $W_{D_n}$) the finite Weyl group of type $A_{n-1}$ (respectively, $B_n$ and $D_n$). Then the group $W_{D_{n}}$ is generated by $s_1,\ldots,s_n$, subject to the following relations:
\begin{align}
s_{i}^{2} & =1\quad (i\leq n-1) \label{eq:invol} \\
s_{i}s_{i+1}s_{i} & =s_{i+1}s_{i}s_{i+1}\quad (i\leq n-2) \\
s_{i}s_{j} & =s_{j}s_{i}\quad (|i-j|>1,\; i,j\neq n) \label{eq:comm} \\
s_{i}s_{n} & =s_{n}s_{i}\quad (i\neq n-2) \\
s_{n-2}s_{n}s_{n-2} & =s_{n}s_{n-2}s_{n}, \quad s_n^2=1.
\label{eq:braidD}
\end{align}
Recall $S_n=W_{A_{n-1}}$ is generated by $s_1,\ldots,s_{n-1}$ subject to the relations (\ref{eq:invol}--\ref{eq:comm}) above.

The group $W_{B_n}$ is generated by $s_1,\ldots,s_n$, subject to the defining relation for $S_n$ on $s_1, \ldots, s_{n-1}$ and the following additional relations:
\begin{align}
s_i s_{n} & =s_ns_i \quad (1\le i \leq n-2) \label{eq:sisnB} \\
(s_{n-1}s_{n})^4 & =1, \; s_{n}^{2} =1. \label{braidB}
\end{align}

\subsection{The extended affine Weyl group $W^e$} In this subsection, we recall the definitions of the extended affine Weyl groups of classical type. For a more detailed expository, consult \cite{Kir}. Define the set $O_W^*$ by
\[
    O_W^* =
    \begin{cases}
        \{1 \},  &\text{if } W=W_{A_{n-1}} \text{ or } W=W_{B_n}\\
        \{n\},  &\text{if } W= W_{D_{2k+1}}\\
        \{1,n\}, &\text{if } W= W_{D_{2k}}.
    \end{cases}
\]

The extended affine Weyl group $W^e$ is the group generated by $s_i$'s and $\pi_r^{\pm 1}$ for $r\in O_W^*$. The generators $s_i$'s obey the relations in (\ref{eq:weyl}) and the defining relations involving $\pi_r$ are shown below:

if $W = W_{A_{n-1}}$, then
\begin{align} \label{eqA:extended weyl}
    \pi_1^2 s_{n-1} &= s_1 \pi_1^2,\nonumber\\
    \pi_1^{n} s_i &= s_{i} \pi_1^{n} \quad (1\le i \le n-1),\\
    \pi_1 s_i &= s_{i+1}\pi_1 \quad (1 \le i \le n-2);\nonumber
\end{align}

if $W = W_{B_n}$, then
\begin{align} \label{eqB:extended weyl}
    \pi_1^2 &= 1, \nonumber\\
    \pi_1 s_i &= s_{i} \pi_1 \quad (2\le i \le n), \\
    \pi_1 s_1 \pi_1 s_1 &= s_1 \pi_1 s_1 \pi_1;\nonumber
\end{align}

if $W=W_{D_n}$ and $n$ is odd, we have
\begin{align} \label{eqDO:extended weyl}
    \pi_n^4 &= 1,  \nonumber\\
    \pi_n^2 s_{n-1} &= s_n \pi_n^2, \\
    \pi_n s_i &= s_{n-i} \pi_n \quad (1\le i \le n-2),\nonumber\\
    \pi_n s_n &= s_1 \pi_n;\nonumber
\end{align}

if $W=W_{D_n}$ and $n$ is even, we have
\begin{align} \label{eqDE:extended weyl}
    \pi_1^2 &= 1, \quad \pi_n^2 = 1, \nonumber\\
    \pi_1 \pi_n &= \pi_n \pi_1, \quad \pi_1 s_1 \pi_1 = \pi_n s_n \pi_n, \nonumber\\
    \pi_1 s_i &= s_i \pi_1, \quad \pi_n s_i = s_{n-i} \pi_n \quad (2\le i \le n-2),\\
    \pi_1 s_{n-1}&=s_n \pi_1, \quad \pi_n s_1 = s_{n-1}\pi_n.\nonumber
\end{align}

\begin{remark}
    We may define
    \[
        s_0 :=
        \begin{cases}
            \pi_1 s_{n-1} \pi_1^{-1}, &\text{for } W= W_{A_{n-1}}\\
            \pi_1 s_{1} \pi_1^{-1}, &\text{for } W= W_{B_{n}} \text{ or } W= W_{D_{2k}}\\
            \pi_n s_{n-1} \pi_n^{-1}, &\text{for } W=W_{D_{2k+1}}.
        \end{cases}
    \] Then the elements $s_i$'s for $i \geq 0$ generate the affine Weyl group $W^a$.
\end{remark}

\subsection{A covering of an extended affine Weyl group}\label{subsec:SpinExtWeyl}
We note here that the Schur multipliers for finite Weyl groups $W$ (and actually for all finite Coxeter groups) have been computed by Ihara and
Yokonuma \cite{IY} (also cf. \cite{Kar}). The explicit generators
and relations for the corresponding covering groups of $W$ can be
found in Karpilovsky \cite[Table 7.1]{Kar}.

In particular, there is a distinguished double covering
$\wtd{W}$ of $W$:
$$
1 \longrightarrow \Z_2 \longrightarrow \wtd{W} \longrightarrow W
\longrightarrow 1.
$$
We denote by $\Z_2 =\{1,z\},$ and by $\td{t}_i$ a fixed preimage
of the generators $s_i$ of $W$ for each $i$. The group $\wtd{W}$
is generated by $z, \td{t}_i$'s with relations
\begin{equation} \label{covering:relations}
z^2 =1, \qquad
 (\td{t}_{i}\td{t}_{j})^{m_{ij}} =
 \left\{
\begin{array}{rl}
1, & \text{if } m_{ij}=1,3 \\
z, & \text{if } m_{ij}=2,4.
\end{array}
\right.
\end{equation}

Let $W$ be $W_{A_{n-1}}$, $W_{B_n}$ or $W_{D_n}$. We define a \textit{covering of the extended affine Weyl group} $W^e$, denoted by $\wtd{W}^e$, to be the group generated by $z, \td{t}_i$'s, and $\td{\pi}_r^{\pm 1}$ ($r \in O_W^*$) such that $z$ is central, $z, \td{t}_i$'s satisfy (\ref{covering:relations}), and the following additional relation:

if $W = W_{A_{n-1}}$, set
\begin{align} \label{eqA:Covering extended weyl}
    \td{\pi}_1^2 \td{t}_{n-1} &= \td{t}_1 \pi_1^2, \nonumber\\
    \td{\pi}_1^{n} \td{t}_i &= \td{t}_{i} \td{\pi}_1^{n}, \\
    \td{\pi}_1 \td{t}_i &= z^{n-1} \td{t}_{i+1}\td{\pi}_1 \quad (1 \le i \le n-2);\nonumber
\end{align}

if $W = W_{B_n}$, set
\begin{align} \label{eqB:Covering extended weyl}
    \td{\pi}_1^2 &= 1, \nonumber\\
    \td{\pi}_1 \td{t}_i &= z \td{t}_{i} \td{\pi}_1 \quad (2\le i \le n), \\
    \td{\pi} \td{t}_1 \td{\pi}_1 \td{t}_1 &= z \td{t}_1 \td{\pi}_1 \td{t}_1 \td{\pi}_1;\nonumber
\end{align}

if $W=W_{D_n}$ and $n$ is odd, set
\begin{align} \label{eqDO:Covering extended weyl}
    \td{\pi}_n^4 &= z, \nonumber\\
    \td{\pi}_n^2 \td{t}_{n-1} &= \td{t}_n \td{\pi}_n^2,\\
    \td{\pi}_n \td{t}_i &=z^{\frac{n-1}{2}} \td{t}_{n-i} \td{\pi}_n \quad (1\le i \le n-2), \nonumber\\
    \td{\pi}_n \td{t}_n &= z^{\frac{n-1}{2}} \td{t}_1 \td{\pi}_n;\nonumber
\end{align}

if $W=W_{D_n}$ and $n$ is even, set
\begin{align} \label{eqDE:Covering extended weyl}
    \td{\pi}_1^2 &= 1, \quad \td{\pi}_n^2 = z^{\frac{n}{2}+1},\nonumber\\
    \td{\pi}_1 \td{\pi}_n &= z \td{\pi}_n \td{\pi}_1, \quad
    \td{\pi}_1 \td{t}_1 \td{\pi}_1 = z\td{\pi}_n \td{t}_n \td{\pi}_n, \nonumber\\
    \td{\pi}_1 \td{t}_i &= \td{t}_i \td{\pi}_1, \quad \td{\pi}_n \td{t}_i = z^{\frac{n}{2}}\td{t}_{n-i} \td{\pi}_n \quad (2\le i \le n-2),\\
    \td{\pi}_1 \td{t}_{n-1}&=\td{t}_n \td{\pi}_1, \quad \td{\pi}_n \td{t}_1 = z^{\frac{n}{2}}\td{t}_{n-1}\td{\pi}_n.\nonumber
\end{align}

The quotient algebra $\C W^{e-} := \C \wtd{W}^e /\langle z+1\rangle$ of
$\C \wtd{W}^e$ by the ideal generated by $z+1$ will be called the
{\em spin extended affine Weyl group algebra} associated to $W$. Denote by $t_i$ (respectively, by $t_{\pi_r}^{\pm 1}$) element in $\C W^{e-}$ the image of $\td{t}_i$ (respectively, $\td{\pi}_r^{\pm 1}$). The {\em spin extended affine Weyl group algebra} $\C W^{e-}$ has the following uniform presentation: $\C W^{e-}$ is the algebra generated by $t_i$'s, and $t_{\pi_r}^{\pm 1}$ for $r \in O_W^*$ subject to the relations
\begin{equation}
(t_{i}t_{j})^{m_{ij}} = (-1)^{m_{ij}+1} \equiv
\begin{cases}
1,  & \text{if } m_{ij}=1,3 \\
-1, & \text{if } m_{ij}=2,4
\end{cases}
\end{equation}
with the following additional relations depending on $W$, i.e, if $W = W_{A_{n-1}}$,
\begin{align} \label{eqA:Spin:extended weyl}
    t_{\pi_1}^2 t_{n-1} &= t_1 t_{\pi_1}^2,\nonumber\\
    t_{\pi_1}^{n} t_i &= t_{i} t_{\pi_1}^{n} \quad (1\le i \le n-1),\\
    t_{\pi_1} t_i &= (-1)^{n-1}t_{i+1}t_{\pi_1} \quad (1 \le i \le n-2);\nonumber
\end{align}

if $W = W_{B_n}$,
\begin{align} \label{eqB:Spin:extended weyl}
    t_{\pi_1}^2 &= 1, \nonumber\\
    t_{\pi_1} t_i &= -t_{i} t_{\pi_1} \quad (2\le i \le n), \\
    t_{\pi_1} t_1 t_{\pi_1} t_1 &= -t_1 t_{\pi_1} t_1 t_{\pi_1} ;\nonumber
\end{align}

if $W=W_{D_n}$ and $n$ is odd,
\begin{align} \label{eqDO:Spin:extended weyl}
    t_{\pi_n}^4 &= -1,  \nonumber\\
    t_{\pi_n}^2 t_{n-1} &= t_n t_{\pi_n}^2,\\
    t_{\pi_n} t_i &= (-1)^{\frac{n-1}{2}}t_{i} t_{\pi_n} \quad (1\le i \le n-2),\nonumber\\
    t_{\pi_n} t_n &= (-1)^{\frac{n-1}{2}} t_1 t_{\pi_n};\nonumber
\end{align}

If $W=W_{D_n}$ and $n$ is even,
\begin{align} \label{eqDE:Spin:extended weyl}
    t_{\pi_1}^2 &= 1, \quad t_{\pi_n}^2 = (-1)^{\frac{n}{2}+1}, \nonumber\\
    t_{\pi_1} t_{\pi_n} &= - t_{\pi_n} t_{\pi_1} \quad
    t_{\pi_1} t_1 t_{\pi_1} = -t_{\pi_n} t_n t_{\pi_n},\nonumber\\
    t_{\pi_1} t_i &= t_i t_{\pi_1} \quad (2\le i \le n-2), \\
    t_{\pi_n} t_i &= (-1)^{\frac{n}{2}}t_{n-i} t_{\pi_n} \quad (2\le i \le n-2),\nonumber \\
    t_{\pi_1} t_{n-1}&=t_n t_{\pi_1}, \quad t_{\pi_n} t_1 = (-1)^{\frac{n}{2}} t_{n-1}t_{\pi_n}.\nonumber
\end{align}

The algebra $\C W^{e-}$ has a superalgebra (i.e. $\Z_2$-graded) structure by letting each $t_i$ be odd and $t_{\pi_r}$ be either even or odd depending on $W$ (motivated by Theorem \ref{th:isofinite} below) as follows: $|t_{\pi_r}| \in \Z_2$, the homogenous degree of $t_{\pi_r}$, is set to be
\begin{align}\label{degree:t_pi_1}
    |t_{\pi_1}| &=
    \begin{cases}
        0, &\text{if } W=W_{A_{2k}} \text{ or } W_{D_{2k}}\\
        1, &\text{if } W= W_{A_{2k+1}} \text{ or } W_{B_n}.
    \end{cases}\quad
\end{align}
\begin{align}\label{degree:t_pi_n}
    |t_{\pi_n}| &=
    \begin{cases}
        k \text{ (mod 2)}, &\text{if } W_{D_n}=W_{D_{2k}}\\
        k \text{ (mod 2)}, &\text{if } W_{D_n}=W_{D_{2k+1}}.
    \end{cases}
\end{align}

%
%
%
%
%
\section{The Clifford algebra} \label{Clifford:Extended Weyl}
In this section, we recall the definition of the Clifford algebra $\Cl_n$. We show that the group $W^e$ acts as automorphisms on $\Cl_n$ which leads to Theorem~\ref{th:isofinite} below.

\subsection{The Clifford algebra $\Cl_n$}
Denote by $\mh^* =\C^n$ the standard (respectively, reflection) representation of the Weyl group $W$ of type $A_{n-1}$ (respectively, of type $B_n$ or $D_n$). Let $\{\al_i\}$ be the set of simple roots for $W$. Note that $\mh^*$ carries a $W$-invariant nondegenerate bilinear form $(-,-)$ such that $(\alpha_i, \alpha_j) =-2\cos (\pi /m_{ij})$. It gives rise to an identification $\mh^*\cong \mh$ and also a bilinear form on $\mh$ which will be again denoted by $(-,-)$.

Denote by $\Cl_n$ the Clifford algebra associated to $(\mh^*,
(-,-))$. We shall denote by $\{c_i\}$ the generators in $\Cl_n$
corresponding to a standard orthonormal basis $\{e_i\}$ of $\C^n$
and denote by $\{\be_i\}$ the elements of $\Cl_n$ corresponding to
the simple roots $\{\al_i\}$ normalized with
$$
\be_i^2=1.
$$
More explicitly, $\mathcal{C}_n$ is generated by
$c_{1},\ldots,c_n$ subject to the relations
\begin{align}  \label{clifford}
c_{i}^{2} =1,\qquad c_{i}c_{j} =-c_{j}c_{i}\quad (i\neq j).
\end{align}
For type $A_{n-1}$, we have
$\be_{i}=\frac{1}{\sqrt{2}}(c_{i}-c_{i+1}),1\leq i\leq n-1$. For
type $B_{n}$, we have an additional $\be_{n}=c_{n}$, and for type
$D_n$, $\be_{n}=\frac{1}{\sqrt{2}}(c_{n-1}+c_{n})$.

The action of $W$ on $\mh$ and $\mh^*$ preserves the bilinear form
$(-,-)$ and thus it acts as automorphisms of the algebra $\Cl_n$.
This gives rise to a semi-direct product $\Cl_n \rtimes \C W$.
Moreover, the algebra $\Cl_n \rtimes \C W$ naturally inherits the
superalgebra structure by letting elements in $W$ be even and each
$c_i$ be odd.

\subsection{An action on $\Cl_n$} \label{action:Cl_n}
We introduce the following elements in $W$. They will be used throughout the paper. Let $\sigma_1^{A} \in W_{A_{n-1}}$ be the cyclic permutation $(1 2 \ldots n)$, $\sigma_1^{D}, \sigma_n^{D} \in W_{D_n},$ and  $\sigma_1^{B} \in W_{B_n}$ be such that

\begin{align*}
\sigma_1^{D} \cdot e_i&=
\begin{cases}
        -e_i, &\text{if } i=1,n\\
        \text{ }e_i, &\text{if } i\neq 1, n.
\end{cases} \qquad
\sigma_n^{D} \cdot e_i &=
\begin{cases}
        -e_{n+1-i}, &\text{if } i\neq n\\
        (-1)^{n-1}e_1, &\text{if } i= n.
\end{cases}\\
\sigma_1^{B} \cdot e_i&=
\begin{cases}
        -e_1, &\text{if } i=1\\
        \text{ }e_i, &\text{if } i\neq 1.
\end{cases}
\end{align*}

\begin{remark}
     We write $\sigma_r$, instead of $\sigma_r^{X}$ where $X \in \{A,B,D\}$ when the context is clear.
\end{remark}

\begin{proposition}\label{extendWeyl:Clifford}
The extended affine Weyl group $W^e$ acts on $\Cl_n$ by the following formulas: $w\cdot c = c^{w}$ and $\pi_r \cdot c = c^{\sigma_r}$ where $c \in \Cl_n, w \in W$ and $r \in O_W^*$.
\end{proposition}
\begin{proof}
     We know that $W$ naturally acts on $\Cl_n$. Suppose $W=W_{A_{n-1}}$, then we have the followings: for $c \in \Cl_n$,
    \begin{eqnarray*}
        \pi_1^2 s_{n-1} \cdot c
        &=& \pi_1^2 \cdot c^{s_{n-1}}
        = c^{\sigma_1^2 s_{n-1}}
        = c^{s_1 \sigma_1^2} = s_1 \pi_1^2 \cdot c.
    \end{eqnarray*}
    For $1\le i \le n-1$,
    \begin{eqnarray*}
        \pi_1^n s_{i} \cdot c
        &=& \pi_1^n \cdot c^{s_{i}}
        = c^{\sigma_1^n s_{i}}
        = c^{s_1 \sigma_1^n} = s_i \pi_1^n\cdot c.
    \end{eqnarray*}
    For $1 \le i \le n-2$,
    \begin{eqnarray*}
        \pi_1 s_{i} \cdot c
        &=& \pi_1 \cdot c^{s_{i}}
        = c^{\sigma_1 s_{i}}
        = c^{s_{i+1} \sigma_1} = s_{i+1} \pi_1\cdot c.
    \end{eqnarray*}
    Therefore, $W^e$ acts on $\Cl_n$. For $W = W_{B_n}$ or $W=W_{D_n}$, the computation is similar to type $A_{n-1}$ case, and can be verified easily.
\end{proof}

Proposition \ref{extendWeyl:Clifford} gives rise to a semi-direct product $\Cl_n \rtimes \C W^e$. Moreover, the algebra $\Cl_n \rtimes \C W^e$ naturally inherits the superalgebra structure by letting elements in $W^e$ be even and each $c_i$ be odd.

\subsection{A superalgebra isomorphism}
Given two superalgebras $\mathcal{A}$ and $\mathcal{B}$, we view
the tensor product of superalgebras $\mathcal{A}$ $\otimes$
$\mathcal{B}$ as a superalgebra with multiplication defined by
\begin{equation}
(a\otimes b)(a^{\prime}\otimes b^{\prime})
=(-1)^{|b||a^{\prime}|}(aa^{\prime }\otimes bb^{\prime})\text{ \ \
\ \ \ \ \ }(a,a^{\prime}\in\mathcal{A},\text{
}b,b^{\prime}\in\mathcal{B})
\end{equation}
where $|b|$ denotes the $\Z_2$-degree of $b$, etc. Also, we shall
use short-hand notation $ab$ for $(a\otimes b) \in \mathcal{A}$
$\otimes$ $\mathcal{B}$, $a = a\otimes1$, and $b=1\otimes b$.

\begin{theorem} \label{th:isofinite}
Let $W=W_{A_{n-1}}, W_{B_n},$ or $W_{D_n}$ and $\upsilon_i = \frac{1}{\sqrt{2}}(c_i + c_{n+1-i})$ for $1\le i \le n$. We have an isomorphism of superalgebras:
$$
\Phi: \Cl_n \rtimes \C W^e
\stackrel{\simeq}{\longrightarrow} \Cl_n \otimes \C W^{e-}
$$
which extends the identity map on $\Cl_n$, sends $s_i \mapsto
-\sqrt{-1} \be_i t_i$, and
\[
    \pi_1 \mapsto
    \begin{cases}
            \be_1 \cdots \be_{n-1} t_{\pi_1}, &\text{ for type $A_{n-1}$}\\
            -\sqrt{-1}c_1t_{\pi_1}, &\text{ for type $B_n$}\\
            \sqrt{-1}c_1c_nt_{\pi_1},&\text{ for type $D_{2k}$},
    \end{cases}
\]

\[
    \pi_n \mapsto
    \begin{cases}
            \upsilon_1 \cdots \upsilon_{\frac{n}{2}} t_{\pi_n}, &\text{ for type $D_n=D_{2k}$}\\
            c_1 \upsilon_1 \cdots \upsilon_{\frac{n-1}{2}}c_{\frac{n+1}{2}} t_{\pi_n}, &\text{ for type $D_n=D_{2k+1}$}.
    \end{cases}
\]
The inverse map $\Psi$ is the extension of
the identity map on $\Cl_W$ which sends $ t_i \mapsto \sqrt{-1}
\be_i s_i$ and

\[
    t_{\pi_1} \mapsto
    \begin{cases}
            \be_{n-1} \cdots \be_{1} \pi_1, &\text{ for type $A_{n-1}$}\\
            \sqrt{-1}c_1 \pi_1, &\text{ for type $B_n$}\\
            -\sqrt{-1}c_nc_1 \pi_1,&\text{ for type $D_n=D_{2k}$},
    \end{cases}
\]

\[
    t_{\pi_n} \mapsto
    \begin{cases}
            \upsilon_{\frac{n}{2}} \cdots \upsilon_{1} \pi_n, &\text{ for type $D_n=D_{2k}$}\\
            c_{\frac{n+1}{2}} \upsilon_{\frac{n-1}{2}} \cdots \upsilon_{1} c_1 \pi_n, &\text{ for type $D_n=D_{2k+1}$}.
    \end{cases}
\]
\end{theorem}
\begin{proof}
By \cite[Theorem ~2.4]{KW1}, the map $\Phi$ (respectively, $\Psi$) preserves the relations not involving $\pi_r$ (respectively, $t_{\pi_r}$) for $r \in O_W^*$. Below, we show that $\Phi$ (respectively, $\Psi$) preserves the relations involving $\pi_r$ (respectively, $t_{\pi_r}$).

If $W = W_{B_{n}}$, we see that $\Phi$ preserves (\ref{eqB:extended weyl}) as follows: for $2 \le i \le n$,
\allowbreak{
\begin{align*}
    \Phi(\pi_1^2)& = -c_1 t_{\pi_1}c_1 t_{\pi_1} = t_{\pi_1}^2 = \Phi(1).\\
    \Phi(\pi_1 s_{1}\pi_1 s_{1}) &= c_1 t_{\pi_1} \be_1 t_1 c_1 t_{\pi_1} \be_1 t_1= c_1 \be_1 c_1 \be_1 t_{\pi_1} t_1 t_{\pi_1} t_1\\
    &= - \be_1 c_1 \be_1 c_1 t_1 t_{\pi_1} t_1 t_{\pi_1}
    = \Phi(s_1 \pi_1s_1 \pi_1).\\
    \Phi(\pi_1 s_{i}) &= -c_1 t_{\pi_1}\be_i t_i = -\be_i c_1 t_{\pi_1}t_i\\
    &= -\be_i t_i c_1 t_{\pi_1}
    = \Phi(s_i \pi_1).
\end{align*}
}
The map $\Psi$ preserves (\ref{eqB:Spin:extended weyl}) as follows: for $2 \le i \le n$,
\allowbreak{
\begin{align*}
    \Psi(t_{\pi_1}^2)& = - c_1 \pi_1 c_1 \pi_1 = \pi_1^2 = \Phi(1).\\
    \Psi(t_{\pi_1} t_1 t_{\pi_1} t_1)
    &= c_1 \pi_1 \be_1 s_1 c_1 \pi_1 \be_1 s_1 \\
    &= c_1(c_1+c_2) \pi_1 s_1 c_1(c_1+c_2) \pi_1 s_1
    = \Psi(-t_1 t_{\pi_1} t_1 t_{\pi_1}).\\
    \Psi(t_{\pi_1} t_{i}) &= -c_1 \pi_1 \be_i s_i = \be_i c_1 \pi_1 s_i = \be_i s_i c_1 \pi_1\\
    &= \Phi(-t_i t_{\pi_1}).
\end{align*}
}
Therefore, $\Phi$ and $\Psi$ are  (super)algebra homomorphisms for $W=W_{B_{n}}$. We see that $\Psi$ is an inverse map of $\Phi$ on the generators. Hence, $\Phi$ and $\Psi$ are inverse algebra isomorphisms.
For $W=W_{A_{n-1}}$ and $W=W_{D_n}$, the computation is similar, but lengthly, and hence will be skipped.
\end{proof}
%
%
%
%
%
%
%
%
%
%
\section{The Trigonometric Double Affine Hecke-Clifford algebra} \label{d DaHC}

In this section, we introduce the trigonometric double
affine Hecke-Clifford algebras, $\dDaHC_W$, and then establish their PBW properties. For $W=W_{A_{n-1}}$, we show that $\dDaHC_W$ is isomorphic to $\daHC_{tr}$ defined in \cite{W}.

\subsection{The algebras $\dDaHC_W$ of type $A_{n-1}$}
Recall $\sigma_1^{A} \in W_{A_{n-1}}$ the cyclic permutation $(1 2 \ldots n)$.
\begin{definition}
Let $u\in \C$, and $W = W_{A_{n-1}}$. The {\em trigonometric double affine Hecke-Clifford algebra}, denoted by $\dDaHC_{W}$ or $\dDaHC_{A_{n-1}}$, is the algebra generated by $\pi_1^{\pm 1},s_1,\ldots,s_{n-1}$, $x_i, c_i$, $(1\le i \le n)$, subject to the relations (\ref{eq:invol}--\ref{eq:comm}), (\ref{clifford}), and the following additional relations:
\allowbreak{
\begin{align}
    x_i x_j & =x_j x_i \quad (\forall i,j) \label{polyn} \\
    x_{i}c_{i} & =-c_{i}x_{i},\text{ }x_{i}c_{j}=c_{j}x_{i}\quad (i\neq j) \label{xici} \\
    \sigma c_{i} & =c_{\sigma i}\sigma\quad (1\leq i\leq
    n,\text{ }\sigma \in S_{n}) \label{sigmac} \\
    x_{i+1}s_{i}-s_{i}x_{i} & =u(1-c_{i+1}c_{i}) \label{xisi} \\
    x_{j}s_{i} & =s_{i}x_{j}\quad (j\neq i,i+1) \label{xjsi}\\
    \pi_1^2 s_{n-1}&= s_1 \pi_1^2 \label{pi^2}\\
    \pi_1 s_i &= s_{i+1} \pi_1 \quad (1\le i \le n-2) \label{pi:s_i}\\
    \pi_1^n s_i &= s_i \pi_1^n \quad (1\le i \le n-1) \label{pi^n:s_i} \\
    \pi_1 x_i &= x_i^{\sigma_1^{A}} \pi_1, \quad \pi_1 c_i = c_i^{\sigma_1^{A}} \pi_1. \label{pi:x_i c_i}
\end{align}
}
\end{definition}
\begin{remark}
    Without $\pi_1^{\pm 1}$, we arrive at the defining relations of the {\em degenerate affine Hecke-Clifford algebra} $\aHC_{A_{n-1}}$ (see \cite{Naz}).
\end{remark}

\subsection{PBW basis for $\dDaHC_{A_{n-1}}$} \label{subsec:PBW}
Let us introduce the ring $\C[P]= \C[P_1^{\pm 1}, \ldots, P_n^{\pm 1}]$. The group $W=S_n$ acts on the ring by the formula
\[
  P_i^{w} = P_{w(i)} \quad \forall w\in S_n.
\]
Also recall the triangular decomposition of $\aHC_W$ in \cite[~Th. 3.4]{KW1}, namely $\aHC_W \cong \C[\mh^*]\otimes \Cl_n \otimes \C W$. We set
\[
    \text{E} := \text{Ind}_{W}^{\aHC_W} \C[P] \cong \C[\mh^*]\otimes \Cl_n \otimes \C[P].
\]
Hence we obtain the representation $\Phi^{\textit{aff}}: \aHC_W \rightarrow \text{End(E)}$ where $x_i$ and $c_i$ act by left multiplication. Together with \cite[Prop. 3.3]{KW1}, we have the followings: for $f \in \C[\mh^*]$, $c\in \Cl_n$, and $g \in \C[P]$, the generator $s_i \in \aHC_W$ for $1\le i \le n-1$ acts by
\begin{align} \label{induced action}
    s_i \cdot (f c g) &= f^{s_i} c^{s_i} g^{s_i} + \left (
    u\disp \frac{f - f^{s_i}}{x_{i+1} - x_i} + u\disp
    \frac{c_ic_{i + 1}f - f^{s_i} c_ic_{i + 1}}{x_{i+1} + x_i} \right
    ) c g.
\end{align}

\begin{lemma}\label{repn:dDaHCa}
    Let $\sigma_1 :=\sigma_1^{A}$. The representation $\Phi^{\textit{aff}}: \aHC_{A_{n-1}} \rightarrow \text{End(E)}$ extends to the representation $\Phi: \dDaHC_{A_{n-1}} \rightarrow \text{End(E)}$ by letting $\pi_1 \cdot (f c g) = P_1 f^{\sigma_1} c^{\sigma_1} g^{\sigma_1}$ where $f \in \C[\mh^*]$, $c\in \Cl_n$, $g \in \C[P]$.
\end{lemma}
\begin{proof}
    We need to check that the relations in $\dDaHC_W$ are preserved under the map $\Phi$. It is enough to check the relations involving only $\pi_1$. Suppose $f \in \C[\mh^*]$, $c\in \Cl_n$, $g \in \C[P]$, then the map preserves the relation (\ref{pi^2}) as follows:
    \allowbreak{
    \begin{eqnarray*}
        \pi_1^2 s_{n-1} \cdot (f c g)
        &=& \pi_1^2 \cdot\left(f^{s_{n-1}} c^{s_{n-1}} g^{s_{n-1}}\right)\\
        && +\pi_1^2 \left( u\disp \frac{f - f^{s_{n-1}}}{x_{n} - x_{n-1}} + u\disp\frac{c_{n-1}c_{n}f - f^{s_{n-1}} c_{n-1}c_{n}}{x_{n} + x_{n-1}} \right) c g\\
        &=& P_1 P_2\left(f^{s_1 \sigma_1^2} c^{s_1 \sigma_1^2} g^{s_1 \sigma_1^2}\right)\\
        && + P_1 P_2 \left( u\disp \frac{f^{\sigma_1^2} - f^{s_{1}\sigma_1^2}}{x_{2} - x_{1}} + u\disp\frac{c_{1}c_{2}f^{\sigma_1^2} - f^{s_{1}\sigma_1^2} c_{1}c_{2}}{x_{2} + x_{1}} \right) c^{\sigma_1^2} g^{\sigma_1^2}\\
        &=& s_1 \pi_1^2 \cdot(f c g).
    \end{eqnarray*}
    } 
    The map preserves the relation (\ref{pi:s_i}) as follows: for $1\le i \le n-2$,
    \begin{eqnarray*}
    \pi_1 s_i \cdot(f c g)
    &=& \pi_1\cdot\left(f^{s_i} c^{s_i} g^{s_i} + \left ( u\disp \frac{f - f^{s_i}}{x_{i+1} - x_i} + u\disp
    \frac{c_ic_{i + 1}f - f^{s_i} c_ic_{i + 1}}{x_{i+1} + x_i} \right
    ) c g \right)\\
    &=& P_1\left(f^{\sigma_1 s_i} c^{\sigma_1 s_i} g^{\sigma_1s_i}\right)\\
    && + P_1\left (u\disp \frac{f^{\sigma_1} - f^{\sigma_1 s_i}}{x_{i+2} - x_{i+1}} + u\disp\frac{c_{i+1}c_{i + 2}f^{\sigma_1} - f^{\sigma_1 s_i} c_{i+1}c_{i + 2}}{x_{i+2} + x_{i+1}} \right) c^{\sigma_1} g^{\sigma_1}\\
    &=& s_{i+1}\pi_1 \cdot (f c g).
    \end{eqnarray*}

    The map preserves the relation (\ref{pi^n:s_i}) as follows: for $1\le i \le n-1$,

    \begin{eqnarray*}
    \pi_1^n s_i \cdot(f c g)
    &=& \pi_1^n\cdot\left(f^{s_i} c^{s_i} g^{s_i} + \left ( u\disp \frac{f - f^{s_i}}{x_{i+1} - x_i} + u\disp
    \frac{c_ic_{i + 1}f - f^{s_i} c_ic_{i + 1}}{x_{i+1} + x_i} \right
    ) c g \right)\\
    &=& P_1P_2\ldots P_n \left(f^{s_i} c^{s_i} g^{s_i}\right)\\
     && + P_1P_2\ldots P_n \left ( u\disp \frac{f - f^{s_i}}{x_{i+1} - x_i} + u\disp\frac{c_ic_{i + 1}f - f^{s_i} c_ic_{i + 1}}{x_{i+1} + x_i} \right) c g\\
    &=& s_i \pi_1^n \cdot(f c g).
    \end{eqnarray*}

    The relation (\ref{pi:x_i c_i}) are can be easily checked.
\end{proof}

\begin{lemma}\label{Faithful S_n^e action}
    Let $W = W_{A_{n-1}}$ and $\sigma_1 :=\sigma_1^{A}$. The extended affine Weyl group $W^e$ acts faithfully on the algebra $\C[\mh^*] \otimes \C[P]$ where $s_i$ for $1\le i \le n-1$ acts naturally and diagonally, while $\pi_1$ acts  as $\sigma_1 \otimes P_1 \sigma_1$.
\end{lemma}
\begin{proof}
    The group $W$ acts diagonally on $\C[\mh^*] \otimes \C[P]$. As operators on $\C[\mh^*] \otimes \C[P]$, we have $\pi_1 s_i \equiv s_{i+1} \pi_1$ for $1\le i \le n-2$, $\pi_1^n s_i \equiv s_i \pi_1^n$ for $1\le i \le n-1$ and $\pi_1^2 s_{n-1} \equiv s_1 \pi_1^2$. Then $W^e$ acts on $\C[\mh^*] \otimes \C[P]$.

    Suppose $\td{w} = \pi_1^{k}s_{i_l}\ldots s_{i_1}$ is a reduced expression with $0 \le i_j \le n-1$ where $s_0 = \pi_1 s_{n-1} \pi_1^{-1}$, and $\td{w}$ acts trivially on $\C[\mh^*] \otimes \C[P]$. Then $\td{w}\cdot 1 = 1$. This implies that $\td{w} \in W$. Since $W$ acts faithfully on $\C[\mh^*] \otimes \C[P]$ this implies that $\td{w} = 1$.
\end{proof}

Below we give a proof of the PBW basis theorem for $\dDaHC_{W}$
using in effect the extension of the representation $\aHC_{W} \rightarrow \text{End(E)}$ to the representation $\dDaHC_{W} \rightarrow \text{End(E)}$. The argument here will be adapted in the proof of Theorem \ref{PBW:DB}.

\begin{theorem} \label{PBW:A}
Let $W=W_{A_{n-1}}$. The set $\{ x^\al c^\be w| \al \in \Z_{+}^n, \beta \in \Z_2^n, w \in W^e\}$ forms a $\C$-linear basis for the $\dDaHC_{W}$. Equivalently, the multiplication of subalgebras
$\C[\mh^*], \Cl_n$, and $\C W^e$ induces a vector space isomorphism
\[
\C[\mh^*]\otimes \Cl_n \otimes \C
W^e\stackrel{\simeq}{\longrightarrow} \dDaHC_W.
\]
\end{theorem}
\begin{proof}
We show that the elements $x^{\al} c^{\be} w$ viewed as operators on $~\C[\mh^*] \otimes \Cl_n \otimes \C[P]$ are linearly independent.

For $\al =(a_1,\ldots,a_n)$, $\nu = (b_1,\ldots,b_n)$, we denote $|\al|=a_1+\cdots+a_n$ and $|\nu|=b_1+\cdots+b_n$. Define a Lexicographic ordering $<$ on the monomials $x^\al P^\nu, \al \in \Z_{\geq 0}^n, \nu \in \Z^n$, by declaring $x^\al P^\nu < x^{\al'} P^{\nu'}$, if $|\al|+|\nu| <|\al'|+|\nu'|$, or if $|\al|+|\nu| =|\al'|+|\nu'|$ then
there exists an $1\le i \le 2n$ such that $z_i<z_i'$ and $z_j=z_j'$
for each $j<i$, where $z_i = a_i$ if $1\le i \le n$, and $z_i = b_{i-n}$ if $n+1 \le i \le 2n$.

Suppose that $S := \sum a_{_{\al \beta w}} x^\al c^\beta w
=0$ for a finite sum over $\al, \beta, w$ and that some
coefficient $a_{_{\al \beta w}} \neq 0$; we fix one such
$\beta$. Now consider the action $S$ on an element of the form
$x_1^{N_1} x_2^{N_2}\cdots x_n^{N_n}$ for $N_n \gg\cdots
\gg N_1 \gg 0$. By Lemma \ref{Faithful S_n^e action}, there exists a unique $\tilde{w}$ such that the leading term of $\td{w}\cdot x_1^{N_1} x_2^{N_2}\cdots x_n^{N_n}$, namely $(x_1^{N_1} x_2^{N_2}\cdots x_n^{N_n})^{\tilde{w}}$, is maximal among all
possible $w$ with $a_{_{\al \beta w}} \neq 0$ for some $\al$. We note that $(x_1^{N_1} x_2^{N_2}\cdots x_n^{N_n})^{\tilde{w}} = (x_1^{N_1} x_2^{N_2}\cdots x_n^{N_n})^{\sigma}P^{\lambda}$ for some $\sigma \in W$ and $\lambda \in \Z^n$. Let $\tilde{\al}$ be chosen among all $\al$ with $ a_{_{\al \be \td{w}}} \neq 0$ such that the monomial ~$x^{\tilde{\al}} (x_1^{N_1} x_2^{N_2}\cdots x_n^{N_n})^{\tilde{w}}c^{\be}$ appears as a maximal term with coefficient $\pm a_{_{{\tilde{\al} \be \tilde{w}}}}$. It follows from $S=0$ that $a_{_{{\tilde{\al} \be \tilde{w}}}} =0$. This is
a contradiction, 
and hence the elements $x^\al c^\be w$ are linearly independent.
\end{proof}
%
%
%
%
%
%
In the classical theory of the trigonometric double affine Hecke algebra  (tDaHa), the algebras are equipped with two presentations. Originally, one defines the tDaHa to be the algebra generated by $\mh^*$ and $W^e$ with certain relations. To obtain the second presentation, one simply rewrite the generators of tDaHa in terms of $\mh^*$, $W$, and the weight lattice $Y$  corresponding to $W$. As an analogue to the classical theory, we show below that the algebra $\daHC_{tr}$ introduced in \cite{W} is isomorphic to $\dDaHC_{A_{n-1}}$.

First, we recall the algebra $\daHC_{tr}$ is generated by $\C[\mh^*]$, $\Cl_n$, $S_n$, and $e^{\pm \epsilon_i}$ $(1 \le i \le n)$, subject to the relations (\ref{xici}--\ref{xjsi}) and the following additional relations:
\begin{align*}
 e^{\epsilon_i} e^{\epsilon_j} &= e^{\epsilon_j} e^{\epsilon_i}, \quad e^{\epsilon_i} e^{-\epsilon_i} =1\\
 w e^{\epsilon_i} &=e^{w\cdot\epsilon_i} w, \quad
 c_j e^{\epsilon_i} =e^{\epsilon_i} c_j  \quad (\forall w\in S_n, \forall i, j) \\
 {[}x_i, e^{\eta}]
 &= u \sum_{k \neq i}  {\text{sgn}(k-i)}
  \frac{e^\eta -e^{s_{ki} (\eta)}}{1 -e^{\text{sgn}(k-i) \cdot (\epsilon_k -\epsilon_i)}} (1 -c_i c_k) s_{ki}.
 \end{align*}

The algebra $\daHC_{tr}$ admits a natural superalgebra structure
with $c_i$ being odd and all other generators being even. Note that the subalgebra generated by $e^{\epsilon_i}$ $(1 \le i \le n)$, denoted by $\C [Y]$ is identified with the group algebra of the weight lattice of type $GL_n$; the subalgebra generated by $e^{\epsilon_i}$ $(1 \le i \le n)$ and $S_n$ is identified with the group algebra of the extended affine Weyl group $W^e$ of type $GL_n$.

\begin{theorem} \label{th:isom dDaHC_A}
Let $W=W_{A_{n-1}}$. We have an isomorphism of superalgebras:
$$
\digamma: \daHC_{tr} \stackrel{\simeq}{\longrightarrow} \dDaHC_W
$$which is the identity map on $\mh^*$, $S_n$, $\Cl_n$, and sends $e^{\epsilon_i} \mapsto s_{i-1}\ldots s_1 \pi_1 s_{n-1} \ldots s_i$ for $1\le i \le n$.
\end{theorem}
\begin{proof}
    By Theorem \ref{PBW:A} and \cite[Prop. 5.4]{W}, it follows that $\daHC_{tr} \cong \dDaHC_W$ as vector spaces. By a long and tedious calculation, we show that the defining relations in $\daHC_{tr}$ are preserved under the map $\digamma$. Hence $\digamma$ is an isomorphism of superalgebras.
\end{proof}

\begin{remark}
    As a consequence of Theorem \ref{PBW:A} and Theorem \ref{th:isom dDaHC_A}, the even center $\mathcal Z(\dDaHC_W)$ contains $\C[Y]^W$ and $\C [x_1^2,\ldots, x_n^2]^W$ as subalgebras. In particular, $\dDaHC_W$ is module-finite over its even center. Alternately, this can be seen for $\daHC_{tr}$ using the locally isomorphism relating to its rational counterpart.
\end{remark}
\subsection{The algebras $\dDaHC_W$ of type $D_{n}$}

Let $W = W_{D_n}$. The extended affine Weyl group $W^e$ is characterized by $n$ even or odd. So the definition of the trigonometric double affine Hecke-Clifford algebra $\dDaHC_W$ of type $D_n$ depends on either $n \in \N$ is even or odd.

Recall that the elements $\sigma_1^D, \sigma_n^D \in W$ are defined in Subsection \ref{action:Cl_n}.

\begin{definition} \label{dDaHCa:D n even}
Let $n \in \N$ be even, $u\in \C$, and $W=W_{D_n}$. The \textit{trigonometric double affine Hecke-Clifford algebra} of type $D_{n}$, denoted by $\dDaHC_W$ or $\dDaHC_{D_n}$, is the algebra generated by and $\pi_1^{\pm 1}, \pi_n^{\pm 1}, x_i, c_i, s_i$, $(1\le i\le n)$ subject to the relations (\ref{eq:invol}--\ref{eq:braidD}), (\ref{polyn}--\ref{xjsi}), and additional relations:
\begin{align}
    s_{n}c_{n} & =-c_{n-1}s_{n}, \quad
    s_{n}c_{i} =c_{i}s_{n} \quad (i\neq n-1,n)\label{sncn}\\
    s_{n}x_{n} & = -x_{n-1}s_{n}-u(1+c_{n-1}c_{n})\\
    s_n x_i & = x_i s_n \quad (i\neq n-1,n)\label{snxi}\\
    \pi_1^2 &=1, \quad \pi_n^2 = 1 \label{pi_1^2=1=pi_n^2}\\
    \pi_1 \pi_n &= \pi_n \pi_1, \quad  \pi_1 s_1 \pi_1 = \pi_n s_n \pi_n \label{pi_1 pi_n=pi_n pi_1}\\
    \pi_1 s_i & = s_i \pi_1, \quad \pi_n s_i = s_{n-i} \pi_n \quad (2\le i \le n-2) \label{pi_1 s_i=s_i pi_1}\\
    \pi_1 s_{n-1} &= s_n \pi_1, \quad \pi_n s_1 = s_{n-1} \pi_n \label{pi_1 s_{n-1}=s_n pi_1}\\
    \pi_r x_i &= x_i^{\sigma_r^D} \pi_r \quad \pi_r c_i = c_i^{\sigma_r^D} \pi_r \quad (r = 1, n; 1\le i \le n).
\end{align}
\end{definition}

\begin{definition} \label{dDaHCa:D n odd}
Let $n \in \N$ be odd, $u\in \C$, and $W=W_{D_n}$. The \textit{trigonometric double affine Hecke-Clifford algebra} of type $D_{n}$, denoted by $\dDaHC_W$ or $\dDaHC_{D_n}$, is the algebra generated by $\pi_n^{\pm 1}, x_i, c_i, s_i$, $(1\le i\le n)$ subject to the relations (\ref{eq:invol}--\ref{eq:braidD}), (\ref{polyn}--\ref{xjsi}), (\ref{sncn}--\ref{snxi}) and additional relations:
\begin{align*}
    \pi_n^4 &= 1  \\
    \pi_n s_i &= s_{n-i} \pi_n \quad (2\le i \le n-2)\\
    \pi_n^2 s_{n-1} & = s_n \pi_n^2 \\
    \pi_n s_1 &= s_{n-1} \pi_n, \quad \pi_n s_n = s_1 \pi_n \\
    \pi_n x_i &= x_i^{\sigma_n^D} \pi_n, \quad \pi_n c_i = c_i^{\sigma_n^D} \pi_n \quad (1\le i \le n).
\end{align*}
\end{definition}

Recall that the element $\sigma_1^B \in W_{B_n}$ is defined in Subsection \ref{action:Cl_n}.

\subsection{The algebras $\dDaHC_W$ of type $B_{n}$}
\begin{definition}
Let $u,v\in \C$, and let $W=W_{B_n}$. The \textit{trigonometric double affine
Hecke-Clifford algebra} of type $B_{n}$, denoted by $\dDaHC_W$ or
$\dDaHC_{B_n}$, is the algebra generated by and $\pi_1^{\pm 1}, x_i, c_i, s_i$, $(1\le i\le n)$ subject to the relations  (\ref{eq:invol}--\ref{eq:comm}), (\ref{eq:sisnB}--\ref{braidB}), (\ref{polyn}--\ref{xjsi}), and additional relations:
\allowbreak{
\begin{align}
    s_{n}c_{n} & =-c_{n}s_{n}, \quad
    s_{n}c_{i} =c_{i}s_{n}\quad (i\neq n) \\
    s_{n}x_{n} & =-x_{n}s_{n} -\sqrt{2} v \\
    s_n x_i &=x_i s_n \quad (i\neq n)\\
    \pi_1^2 &= 1  \label{pi^2:B}\\
    \pi_1 s_i &= s_{i} \pi_1 \quad (2\le i \le n) \label{pi:s_i:B}\\
    \pi_1 s_1 \pi_1 s_1 &= s_1 \pi_1 s_1 \pi_1 \label{pi:s_1} \\
    \pi_1 x_i &= x_i^{\sigma_1^B} \pi_1, \quad \pi_1 c_i = c_{i}^{\sigma_1^B} \pi_1 \quad (1\le i \le n)\label{piB:x_i}.
\end{align}
}
\end{definition}

\begin{remark}
    Without $\pi_r^{\pm 1}$, we arrive at the defining relations of the {\em degenerate affine Hecke-Clifford algebra} $\aHC_W$, for $W=W_{B_{n}}$ or $W_{D_n}$ (see \cite{KW1}).
\end{remark}
\subsection{PBW basis for $\dDaHC_W$}
In this subsection, we prove a PBW type result for the algebra $\dDaHC_W$ for $W=W_{B_n}$ or $W=W_{D_n}$. We first make a suitable modification from Subsection \ref{subsec:PBW}. We introduce the ring $\C[P]:= \C[P_1^{\pm 1/2},\ldots, P_n^{\pm 1/2}]$ where $P_i$ formally corresponds to $e^{\epsilon_i}$. $S_n$ naturally acts on $\C[P]$ by the formula $(P_i^{1/2})^w = P_{w(i)}^{1/2}$. The action can be extended to an action of $W=W_{D_n}$ by letting
\begin{align*}
    (P_n^{1/2})^{s_n} &= P_{n-1}^{-1/2}, \quad (P_{n-1}^{1/2})^{s_n} = P_{n}^{-1/2}\\
    (P_i^{1/2})^{s_n} &= P_i^{1/2}, \quad (i\neq n-1,n).
\end{align*}

Also, the action of $S_n$ on $\C[P]$ can be extended to an action of $W=W_{B_n}$ by letting
\begin{align*}
    (P_n^{1/2})^{s_n} &= P_n^{-1/2}, \quad (P_i^{1/2})^{s_n} = P_i^{1/2}, \quad (i\neq n).
\end{align*}

Consider the representation $\Phi^{\textit{aff}}: \aHC_W \rightarrow \text{End(E)}$, namely
\[
    \text{E} := \text{Ind}_{W}^{\aHC_W} \C[P] \cong \C[\mh^*]\otimes \Cl_n \otimes \C[P].
\]
Note that $x_i$ and $c_i$ act by left multiplication. The generators $s_i \in \aHC_W$ for $W = W_{D_{n}}$ $(1\le i \le n)$ act by the same formula (\ref{induced action}) for $1\le i \le n-1$ and in addition by
\begin{align*}
s_{n}\cdot (fcg)=f^{s_n} c^{s_{n}} g^{s_n}
-\left(u\disp\frac{f-f^{s_n}}{x_{n}+x_{n-1}} -
u\displaystyle\frac{c_{n-1}c_{n}f-f^{s_n}c_{n-1}c_{n}}{x_{n}-x_{n-1}}\right )c g.
\end{align*}

The generators $s_i \in \aHC_W$ for $W = W_{B_{n}}$ $(1\le i \le n)$ act by the same formula (\ref{induced action}) for $1\le i \le n-1$ and in addition by
\begin{align*}
s_{n}\cdot (fcg)=f^{s_n} c^{s_{n}} g^{s_n}
-\sqrt{2}\left( \disp v\frac{f-f^{s_n}}{2 x_{n}}\right )c g.
\end{align*}

For more details treatment on $\aHC_W$, consult \cite{Naz,KW1}.
The following lemma is a counterpart to Lemma \ref{repn:dDaHCa}.

\begin{lemma}\label{repnDB:dDaHCa}
    Let $W=W_{D_n}$ or $W_{B_n}$. The representation $~\Phi^{\textit{aff}}:\aHC_{W}\rightarrow\text{End(E)}$ extends to the representation $\Phi:\dDaHC_{W}\rightarrow \text{End(E)}$ by the following formulas:
    \[
            \pi_r \cdot (f c g) =
            \begin{cases}
                P_1 f^{\sigma_1} c^{\sigma_1} g^{\sigma_1}, &\text{ if } r =1,\\
                (P_1\ldots P_n)^{1/2} f^{\sigma_n} c^{\sigma_n} g^{\sigma_n}, &\text{ if } r =n,
            \end{cases}
        \]
where $f \in \C[\mh^*], c\in \Cl_n, g \in \C[P]$.
\end{lemma}
\begin{proof}
By a direct and lengthy computation, the action $\aHC_W$ on $E \cong \C[\mh^*] \otimes \Cl_n \otimes \C[P]$ naturally extends to an action of $\dDaHC_W$. We will verify a few relations in $\dDaHC_{D_n}$ for $n$ is even, and leave the rest to the reader.

The map preserves the relation (\ref{pi_1^2=1=pi_n^2}) as follows:
\begin{align*}
    \pi_1^2\cdot (fcg) &= \pi_1 \cdot P_1 (fcg)^{\sigma_1} = 1\cdot(fcg).\\
    \pi_n^2 \cdot(f c g)
    &= \pi_n \cdot (P_1\ldots P_n)^{1/2} (fcg)^{\sigma_n} \\
    &= (P_1\ldots P_n)^{1/2} (P_n^{-1}\ldots P_1^{-1})^{1/2} f c g
    = 1\cdot(fcg).
\end{align*}

The map preserves the relation (\ref{pi_1 pi_n=pi_n pi_1}) as follows: since $\sigma_1 \sigma_n = \sigma_n \sigma_1$, then $\pi_1 \pi_n \cdot fcg = \pi_n \pi_1 \cdot fcg$. We also have
\begin{eqnarray*}
\pi_1 s_1 \pi_1 \cdot f c g&=& \pi_1 s_1 \cdot P_1 (fcg)^{\sigma_1}\\
&=& \pi_1 \cdot P_2(fcg)^{s_1 \sigma_1} \\
&& + u \pi_1 \cdot P_2 \left(\disp \frac{f^{\sigma_1} - f^{s_1\sigma_1}}{x_{2} - x_1} + \disp \frac{c_1c_{2}f^{\sigma_1} - f^{s_1\sigma_1} c_1c_{2}}{x_{2} + x_1} \right) (c g)^{\sigma_1}\\
&=& P_1 P_2 (fcg)^{\sigma_1 s_1 \sigma_1} \\
&&+ u P_1 P_2 \left( \disp \frac{f - f^{\sigma_1 s_1\sigma_1}}{x_{2} + x_1} + \disp \frac{c_{2}c_1f - f^{\sigma_1s_1\sigma_1}c_{2}c_1}{x_{2} - x_1} \right) c g.
\end{eqnarray*}

$\pi_n s_n \pi_n \cdot f c g $
\begin{eqnarray*}
&=& \pi_n s_n \cdot (P_1\ldots P_n)^{1/2} (f c g)^{\sigma_n}\\
&=& \pi_n \cdot (P_1\ldots P_{n-2} P_{n-1}^{-1} P_n^{-1})^{1/2} (f c g)^{s_n \sigma_n}\\
&& - u \pi_n \cdot(P_1\ldots P_{n-2} P_{n-1}^{-1} P_n^{-1})^{1/2}\left(\disp \frac{f^{\sigma_n} - f^{s_n\sigma_n}}{x_{n} + x_{n-1}}\right)(c g)^{\sigma_n}\\
&& + u \pi_n \cdot(P_1\ldots P_{n-2} P_{n-1}^{-1} P_n^{-1})^{1/2}\left(\disp \frac{c_{n-1}c_{n}f^{\sigma_n} - f^{s_n\sigma_n} c_{n-1}c_{n}}{x_{n} - x_{n-1}} \right) (c g)^{\sigma_n}\\
&=&P_1 P_2 (f c g)^{\sigma_n s_n \sigma_n} + u P_1 P_2 \left(u\disp \frac{f^ - f^{\sigma_n s_n\sigma_n}}{x_{1} + x_{2}}+  \frac{c_{2}c_{1}f - f^{\sigma_n s_n\sigma_n} c_{2}c_{1}}{-x_{1} + x_{2}} \right) c g.
\end{eqnarray*}
Since $\sigma_1 s_1 \sigma_1 = \sigma_n s_n \sigma_n$, then it follows that $\pi_1 s_1 \pi_1 \cdot f c g = \pi_n s_n \pi_n\cdot f c g.$

The map preserves the relation (\ref{pi_1 s_i=s_i pi_1}) as follows: for $2\le i \le n-2$, we have
\allowbreak{
\begin{eqnarray*}
\pi_n s_i \cdot fcg &=& \pi_n \cdot (fcg)^{s_i} + \pi_n \cdot \left (
u\disp \frac{f - f^{s_i}}{x_{i+1} - x_i} + u\disp
\frac{c_ic_{i + 1}f - f^{s_i} c_ic_{i + 1}}{x_{i+1} + x_i} \right
) c g\\
&=& (P_1\ldots P_{n-1} P_n)^{1/2}(fcg)^{\sigma_n s_i} \\
&& +u (P_1\ldots P_{n-1} P_n)^{1/2} \left (\disp \frac{f^{\sigma_n} - f^{\sigma_n s_i}}{-x_{n-i} + x_{n+1-i}}\right) (c g)^{\sigma_n}\\
&& + u (P_1\ldots P_{n-1} P_n)^{1/2} \left(\disp
\frac{c_{n+1-i}c_{n-i}f - f^{\sigma_n s_i} c_{n+1-i}c_{n-i}}{-x_{n-i}  -x_{n+1-i}} \right) (c g)^{\sigma_n}\\
&=& s_{n-i} \pi_n \cdot fcg.
\end{eqnarray*}
}
It is easy to verify $\pi_1 s_i = s_i \pi_1$ for $2\le i \le n-2$. For the rest of the relations in $\dDaHC_W$, the calculation is similar.
\end{proof}
We have the following PBW basis theorem for $\dDaHC_W$.

\begin{theorem} \label{PBW:DB}
Let $W=W_{D_n}$ or $W=W_{B_n}$, the elements $\{x^\al c^\be w| \al\in\Z_+^n, \be \in\Z_2^n, w\in W^e\}$ form a $\C$-linear basis for $\dDaHC_W$ (called a PBW basis). Equivalently, the multiplication of subalgebras $\C[\mh^*], \Cl_n$, and $\C W^e$ induces a vector space isomorphism
\[
\C[\mh^*]\otimes \Cl_n \otimes \C
W^e\stackrel{\simeq}{\longrightarrow} \dDaHC_W.
\]
\end{theorem}
\begin{proof}

We show that the elements $x^{\al} c^{\be} w$ viewed as operators on $~\C[\mh^*] \otimes \Cl_n \otimes \C[P]$ are linearly independent. It is clear that, for either $W=W_{D_n}$ or $W=W_{B_n}$, the elements $x^\al
c^\be w$ span $\dDaHC_W$. It remains to show that they are linearly
independent. We shall treat the $W_{B_n}$ case in detail and skip
the analogous $W_{D_n}$ case.

To that end, we shall refer to the argument in the proof of
Theorem~\ref{PBW:A} with suitable modification. First, we observe that for each $w \in W_{B_n}^e$, the leading term of $w \cdot (x_1^{N_1} \ldots x_n^{N_n})$ is $(x_1^{N_1} \ldots x_n^{N_n})^{\gamma}P^{\lambda}$ for some $\gamma \in W_{B_n}$ and $\lambda \in (\frac{\Z}{2})^n$. Similar to type $A$ case, we choose $\tilde{w}$ such that the leading term $(x_1^{N_1} x_2^{N_2}\cdots x_n^{N_n})^{\tilde{w}} =  (x_1^{N_1} \ldots x_n^{N_n})^{\gamma} P^{\lambda}$ is maximal for some $\gamma, \lambda$. Write $\gamma = ((\eta_1, \dots, \eta_n), \sigma) \in W_{B_n} =\{\pm 1\}^n \rtimes S_n$. Let $\tilde{\al}$ be chosen among all $\al$ with $ a_{_{\al \be \td{w}}} \neq 0$ such that the monomial ~$x^{\tilde{\al}} (x_1^{N_1} x_2^{N_2}\cdots x_n^{N_n})^{\tilde{w}}c^{\be}$ appears as a maximal term with coefficient $\pm a_{_{{\tilde{\al} \be \tilde{w}}}}$. Note that $\tilde{w}$ may now not be unique, but the $\lambda, \sigma,$ and $\tilde{\al}$ are uniquely determined. Then, by the same argument on the vanishing of a maximal term, we obtain that
\allowbreak{
\begin{eqnarray*}
0 &=&  \sum_{\tilde{w}} a_{{\tilde{\al} \be \tilde{w}}} x^{\tilde{\al}} (x_1^{N_1} x_2^{N_2}\cdots x_n^{N_n})^{\tilde{w}}\\
&=& \sum_{\gamma} a_{{\tilde{\al} \be \tilde{w}}}  x^{\tilde{\al}} (x_1^{N_1} x_2^{N_2}\cdots x_n^{N_n})^{\gamma}P^{\lambda},
\end{eqnarray*}
}
and hence
\[
\sum_{(\eta_1, \dots, \eta_n)} a_{{\tilde{\al} \be
\tilde{w}}} (-1)^{\sum_{i=1}^n \eta_i N_i} =0.
\]

By choosing $N_1,\ldots, N_n$ with different parities ($2^n$
choices) and solving the $2^n$ linear equations, we see that all
$a_{{\tilde{\al} \be \tilde{w}}} =0$. This can also be seen
more explicitly by induction on $n$. By choosing $N_n$ to be even
and odd, we deduce that for a fixed $\eta_n$,
$\sum_{(\eta_1, \dots, \eta_{n-1}) \in \{\pm 1\}^{n-1}}
a_{{\tilde{\al} \epsilon \tilde{w}}} (-1)^{\sum_{i=1}^{n-1} \eta_i
N_i} =0$, which is the equation for $(n-1)$ $x_i$'s and the
induction applies.
\end{proof}

\begin{corollary} \label{center:BD}
Let $W=W_{D_n}$ or $W=W_{B_n}$. The even center of $\dDaHC_W$ contains $\C[x_1^2,\ldots, x_n^2]^W$.
\end{corollary}
\begin{proof}
Suppose $f \in \C [x_1^2,\ldots, x_n^2]^W$. Then $\pi_r f = f^{\sigma_r} \pi_1 = f \pi_1$. By \cite[Prop. 4.6]{KW1}, $f$ commutes with the generators $s_1,\ldots, s_{n}$ and $c_1,\ldots, c_n$. So, $f$ is in the even center of $\dDaHC_W)$.
\end{proof}
%
%
%
\section{Trigonometric spin double affine Hecke algebras} \label{sec:spin}

In this section we introduce the trigonometric spin double affine
Hecke algebra $\sDaH_W$, and then establish its connections to the corresponding trigonometric double affine Hecke-Clifford algebras $\dDaHC_W$.

\subsection{The skew-polynomial algebra}
We shall denote by $\Cl[\xi_1,\ldots,\xi_n]$ the $\C$-algebra
generated by $\xi_1,\ldots,\xi_n$ subject to the relations
\begin{align}
\xi_i \xi_j + \xi_j \xi_i =0\quad (i\neq j). \label{skew:polyn}
\end{align}
This is naturally a superalgebra by letting each $\xi_i$ be odd.
We will refer to this as the {\em skew-polynomial algebra} in $n$
variables. This algebra has a linear basis given by $\xi^\alpha
:=\xi_1^{k_1}\cdots \xi_n^{k_n}$ for $\alpha =(k_1,\ldots,k_n) \in
\Z_+^n$, and it contains the polynomial subalgebra
$\C[\xi_1^2,\ldots, \xi_n^2]$.

\subsection{The algebra $\sDaH_W$ of type $A_{n-1}$}
Let $W=W_{A_{n-1}}$.  Recall that the spin extended affine Weyl group $\C W^{e-}$ associated to a Weyl group $W$ is generated by $t_1,\ldots,t_{n-1}$, and $t_{\pi_1}$ subject to the relations as specified in Section ~\ref{subsec:SpinExtWeyl}

\begin{definition}
Let $u \in \C$. The {\em trigonometric spin double affine Hecke algebra} of type $A_{n-1}$, denoted by $\sDaH_W$ or $\sDaH_{A_{n-1}}$, is the algebra generated by $\Cl[\xi_1,\ldots,\xi_n]$ and $\C W^{e-}$ subject to the relations:
\begin{align}
t_{i}\xi_{i} &= - \xi_{i+1}t_i +u  \quad (1\le i\le n-1) \label{ti xi_i:A} \\
t_i\xi_j &= - \xi_{j}t_i \quad (j\neq i,i+1) \label{ti xi_j:A}\\
t_{\pi_1} \xi_i &= (-1)^{n-1}\xi_{i+1} t_{\pi_1} \quad (1\le i \le n-1)\\
t_{\pi_1} \xi_n &= (-1)^{n-1} \xi_1 t_{\pi_1}.
\end{align}
\end{definition}

\subsection{The algebra $\sDaH_W$ of type $D_n$}
Let $u \in \C$ and $W=W_{D_n}$. In this subsection we define the {\em trigonometric spin double affine Hecke algebra} associated to $W$, denoted by $\sDaH_W$ or $\sDaH_{D_n}$, for both $n$ is even or odd.

\begin{definition} \label{def:sDaHD}
Let $n$ be odd. The algebra $\sDaH_W$ is generated by $\Cl[\xi_1,\ldots,\xi_n]$ and $\C W^{e-}$ subject to the relations(\ref{ti xi_i:A}-\ref{ti xi_j:A}) and the following additional relations:
\begin{align}
t_n\xi_n &= - \xi_{n-1}t_n + u\\
t_n\xi_i &=-\xi_i t_n \quad (i \neq n-1,n)\label{t_n xi_i} \\
t_{\pi_n} \xi_i &= (-1)^{\frac{n-1}{2}} \xi_{n+1-i} t_{\pi_n} \quad (1 \le i \le n).
\end{align}
\end{definition}

\begin{definition}
Let $n$ even, the algebra $\sDaH_W$ is generated by $\Cl[\xi_1,\ldots,\xi_n]$ and $\C W^{e-}$ subject to the relations (\ref{ti xi_i:A}-\ref{t_n xi_i}) and the following additional relations:
\begin{align*}
    t_{\pi_1} \xi_i &= \xi_i t_{\pi_1} \quad (1\le i \le n)\\
    t_{\pi_n} \xi_i &= (-1)^{\frac{n}{2}+1}\xi_{n+1-i} t_{\pi_n}.
\end{align*}
\end{definition}

\subsection{The algebra $\sDaH_W$ of type $B_n$}

\begin{definition} Let $u,v\in\C$, and $W=W_{B_n}$.
The {\em trigonometric spin double affine Hecke algebra} of type $B_n$, denoted by $\sDaH_W$ or $\sDaH_{B_n}$, is the algebra generated by
$\Cl[\xi_1,\ldots,\xi_n]$ and $\C W^{e-}$ subject to the relations
relations (\ref{ti xi_i:A}-\ref{ti xi_j:A}) and the following additional relations:
\begin{align*}
t_n\xi_n &= -\xi_{n}t_n + v\\
t_n\xi_i &=-\xi_{i}t_n \quad (i\neq n)\\
t_{\pi_1} \xi_i & = -\xi_i t_{\pi_1} \quad (1\le i\le n).
\end{align*}
\end{definition}

Sometimes, we will write $\sDaH_W(u,v)$ or $\sDaH_{B_n}(u,v)$ for
$\sDaH_W$ or $\sDaH_{B_n}$ to indicate the dependence on the
parameters $u,v$.

\begin{remark}
    The algebra $\sDaH_W$ is naturally a superalgebra with all $t_i$'s and $\xi_i$ being odd generators. A generator $t_{\pi_r}$ can be either even or odd depending on $r$ and $W$, its degree is given by (\ref{degree:t_pi_1}) and (\ref{degree:t_pi_n}).
    Without $t_{\pi_r}$, we arrive at the defining relation of {\em degenerate spin affine Hecke algebra}, denoted by $\saH_W$, see \cite{W, KW1}.
\end{remark}

\subsection{A superalgebra isomorphism}
\begin{theorem} \label{th:isomABD}
Let $W=W_{A_{n-1}}, W_{D_n}$ or $W_{B_n}$. Then,
\begin{enumerate}
\item there exists an isomorphism of superalgebras
$$\Phi:\dDaHC_W {\longrightarrow }\Cl_n \otimes \sDaH_W$$
which extends the isomorphism $\Phi: \Cl_n \rtimes \C W^{e}
\longrightarrow \Cl_n \otimes \C W^{e-}$ (in Theorem~\ref{th:isofinite}) and sends
$x_{i}\longmapsto\displaystyle \sqrt{-2}c_{i}\xi_{i}$ for each
$i;$
\item the inverse $\Psi:\mathcal{C}_{n}\otimes
\sDaH_W{\longrightarrow}\dDaHC_W$ extends $\Psi:\Cl_n\otimes\C
W^{e-} \longrightarrow \Cl_n \rtimes \C W^{e}$ (in
Theorem~\ref{th:isofinite}) and sends $\xi_{i}\longmapsto
\displaystyle\frac{1}{\sqrt{-2}}c_{i}x_{i}$ for each $i$.
\end{enumerate}
\end{theorem}
\begin{proof}
    We need to show that $\Phi$ (respectively, $\Psi$) preserves the defining relations in $\dDaHC_W$ (respectively, in $\sDaH_W$). By \cite[Theorem 4.4]{KW1}, the map $\Phi$ (respectively, $\Psi$) preserves the relations not involving only $\pi_r$ (respectively, $t_{\pi_r}$). So it is left to show that $\Phi$ (respectively, $\Psi$) preserves the defining relations which involve $\pi_r$ and $x_i$'s (respectively, $t_{\pi_r}$ or $\xi_i$'s) for $1\le i \le n$. We verify relations for type $B_n$ case below. For other relations and types, the computation is similar, and will be skipped.

    For $W = W_{B_{n}}$ and $2\le i \le n$, we have
    \begin{align*}
    \Phi(\pi_1 x_i) & = -\sqrt{2} c_1 t_{\pi_1} c_i \xi_i =  -\sqrt{2} c_i c_1 t_{\pi_1} \xi_i = \Phi(x_i \pi_1).\\
    %
    \Phi(\pi_1 x_1)
        & = -\sqrt{2}c_1 t_{\pi_1} c_1 \xi_1 =\sqrt{2} t_{\pi_1}\xi_1
        =  -\sqrt{2} \xi_1 t_{\pi_1} = \sqrt{2} c_1 \xi_1 c_1t_{\pi_1}\\
        & = \Phi(-x_1 \pi_1).
    \end{align*}
Thus $\Phi$ is a homomorphism of (super)algebras. Similarly, $\Psi$ is a superalgebra homomorphism by the followings: for $1\le i \le n$
\begin{align*}
    \Psi(t_{\pi_1} \xi_i)
    & = \frac{1}{\sqrt{2}} c_1 \pi_1 c_i x_i =  -\frac{1}{\sqrt{2}} c_i x_i c_1 \pi_1  = \Psi(-\xi_i t_{\pi_1}).
\end{align*}

Since $\Phi$ and $\Psi$ are inverses on generators and hence they are indeed (inverse) isomorphisms.
\end{proof}

\subsection{PBW basis for $\sDaH_W$}
We have the following PBW basis theorem for $\saH_W$.

\begin{theorem} \label{PBW:spin}
Let $W=W_{A_{n-1}}$, $W_{D_n}$ or $W_{B_n}$. The multiplication of the
subalgebras $\C W^-$ and $\Cl [\xi_1,\ldots,\xi_n]$ induces a vector
space isomorphism
\[
\Cl [\xi_1,\ldots,\xi_n] \otimes \C W^{e-}
\stackrel{\simeq}{\longrightarrow} \sDaH_W.
\]
\end{theorem}

\begin{proof}
It follows from the definition that $\sDaH_W$ is spanned by the
elements of the form $\xi^{\al} \sigma$ where $\sigma$ runs over a
basis for $\C W^{e-}$ and $\al \in \Z_+^n$. By
Theorem~\ref{th:isomABD}, we have an isomorphism
$\Psi:\Cl_{n}\otimes \sDaH_W {\longrightarrow}\dDaHC_W.$ Observe
that the image $\Psi(\xi^{\al} \sigma)$ are linearly independent in
$\dDaHC_W$ by the PBW basis Theorems~\ref{PBW:A} and \ref{PBW:DB} for $\dDaHC_W$. Hence the elements $\xi^{\al} \sigma$ are linearly independent in $\sDaH_W$.
\end{proof}

\begin{remark}
$\sDaH_W$ contains the skew-polynomial algebra $\Cl
[\xi_1,\ldots,\xi_n]$ and the spin extended affine Weyl group algebra $\C W^{e-}$ as subalgebras.
\end{remark}

\begin{remark}
As a counterpart of $\daHC_{tr}$, \cite{W} also introduce the trigonometric spin double affine Hecke-Clifford algebra of type $A$, denoted by $\sdaH_{tr}$. As a consequence of Theorems \ref{th:isom dDaHC_A} and \ref{th:isomABD}, we have $\sDaH_{A_{n-1}} \cong \sdaH_{tr}$, and therefore $\sDaH_{A_{n-1}}$ is module-finite over its even center. We also have the following.
\end{remark}
\begin{corollary}
    Let $W = W_{A_{n-1}}, W_{B_n}$ or $W_{D_n}$. The even center for $\sDaH_W$ contains $\C [\xi_1^2,\ldots, \xi_n^2]^W$.
\end{corollary}
\begin{proof}
    By the isomorphism $\Phi:\dDaHC_W {\longrightarrow }\Cl_n \otimes \sDaH_W$, we have
\[
Z(\Cl_n \otimes \sDaH_W) =\Phi(Z(\dDaHC_W)) \supseteq \Phi(\C [x_1^2,\ldots,x_n^2]^W) =\C [\xi_1^2,\ldots,\xi_n^2]^W.
\]
Thus, $\C [\xi_1^2,\ldots,\xi_n^2]^W \subseteq Z(\sDaH_W)$.
\end{proof}

%

%
%

%
%
%
%
%
%

\end{document}